%% file: main.tex
\documentclass[journal]{new-aiaa} 
\usepackage[utf8]{inputenc}

\usepackage{graphicx}
\usepackage{amsmath}
\usepackage[version=4]{mhchem}
\usepackage{siunitx}
\usepackage{longtable,tabularx}
\usepackage{tikz, pgfplots}
\usepackage{pgfplots, pgfplotstable}
\usepackage{filecontents}
\usepackage{subcaption}
\title{Space Logistics Modeling and Optimization: \\Review of the State of the Art\footnote{The previous version of this paper was presented at the AIAA SciTech Forum in Orlando, FL, on Jan 8-12 (AIAA 2024-1275).}}

\author{Koki Ho\footnote{Dutton-Ducoffe Professor, Associate Professor, Daniel Guggenheim School of Aerospace Engineering, Chair of the AIAA Space Logistics Technical Committee, AIAA Senior Member}}
\affil{Georgia Institute of Technology, Atlanta, GA, 30309}

\begin{document}

\maketitle

\section{Introduction}
The US Space Force's doctrine document published in 2020, \textit{SpacePower: Doctrine for Space Forces} listed five core competencies that the US Space Forces need to perform, one of which was ``Space Mobility and Logistics'' \cite{Raymond2020}. In response to that, there has been a rapidly growing interest in developing technologies to enable in-space operations via in-space infrastructures and in-space refueling, servicing, assembly, and manufacturing (ISAM). Space logistics has also been a topic of interest to NASA and other space agencies as an approach to designing sustainable space exploration campaigns. Additionally, there is a growing need for space logistics to effectively operate the mega-scale constellations in Earth orbits for communication, remote sensing, and related missions.

Enabling space mobility and logistics capabilities requires a new way to view space missions that differs from conventional astrodynamics. For example, with fuel depots in space, the optimal path to the destination in terms of the lifecycle cost is not necessarily a fuel-optimal trajectory; rather, a path stopping by a fuel depot and being refueled before heading to the destination may be preferred even when it requires additional fuel. In addition, with ISAM capabilities, we need additional analysis capabilities to analyze and optimize the sizes of the fuel/spare depots and their inventory/sparing policies with astrodynamics in mind. This analysis requires logistics-driven modeling and optimization techniques coupled with astrodynamics. This intersection between these areas leads to an emerging field, Space Logistics modeling and optimization.  

Space Logistics is formally defined as ``the theory and practice of driving space system design for operability and supportability, and of managing the flow of materiel, services, and information needed throughout a space system lifecycle \cite{AIAASLTC}''  according to the AIAA Space Logistics Technical Committee. Space Logistics is not a new research topic; the literature on space logistics has existed since the beginning of the space exploration era \cite{freeman1966mathematical,rootevaluation,carrillo1983development}. According to the AIAA's meeting papers archive, the 1st AIAA/SOLE Space Logistics Symposium was held in 1987 \cite{AIAASOLE}, and it already had papers that covered various topics that are still relevant today such as in-space servicing \cite{shepard1987role,bell1987logistics} and even using artificial intelligence (AI) for space logistics \cite{barryapplication}. Although the applications of interest back then do not necessarily match with our current ones, various key concepts were developed decades ago. In the early 2000s, in response to NASA's Constellation Program, the idea of logistics network modeling was applied to human space mission design. In addition, more logistics-driven techniques for probabilistic modeling and inventory control have been applied to satellite servicing and mega-scale constellations. A recently published book, \textit{The Planning and Execution of Human Missions to the Moon and Mars} \cite{poliskie2023planning}, has two chapters that focus on space logistics and mathematical modeling for space mission design \cite{shishko2023interplanetary,ho2023mathematical}, which can be referenced for more details.

The goal of this paper is to categorize the state-of-the-art studies in Space Logistics modeling and optimization in two ways: (1) by application questions that are addressed; and (2) by logistics-driven methods that are used in the studies. The goal of the first categorization is to help the practitioners determine what research is out there that can support their applications. The goal of the second categorization is to systematically map the existing literature to each logistics research subfield, and thus help the researchers to understand the state of the art and identify the under-explored and promising future research directions.

\section{Categorization of State of the Art by Applications}
\label{sec:app_categorization}
In this section, we first categorize the state-of-the-art literature by applications to aid the practitioners in understanding the existing research that can answer their application questions. We will introduce the three major applications in Space Logistics, and then map the state of the art in the literature to each research question. 

Note that this is by no means a comprehensive list of all space logistics applications. This paper specifically focuses on the applications that involve orbital mechanics considerations; other space logistics applications not reviewed in this paper include supportability and spare parts management within a space station \cite{siddiqi2007spare,owens2017supportability,kline2007estimating}, surface vehicle routing and logistics \cite{siddiqi2006reconfigurability,ahn2008optimization,ahn2010mars,Ahn2012,leeD2018,leeD2019,choi2022,choi2023}, among others.

\subsection{Application 1 (A1): In-space Servicing, Assembly, and Manufacturing (ISAM) for Satellites}
Although In-space Servicing, Assembly, and Manufacturing (ISAM) has attracted a lot of attention in recent years, the concepts themselves have been explored for a long time. When an operational satellite experiences a component failure or runs out of its fuel, it could be economical to repair and refuel the satellite on-orbit than launch a new one from the ground especially when the satellite is in a high-altitude orbit. One of the earlier successful examples of on-orbit servicing was the repair mission for the Hubble Space Telescope, although crewed repair missions have become difficult after the retirement of the Space Shuttle. In the 2000s, multiple studies introduced the modeling and simulation foundations for robotic/autonomous on-orbit servicing missions \cite{saleh2002space,lamassoure2002space,long2007orbit} and examined the feasibility and economic analysis of various ISAM concepts \cite{galabova2006economic,gralla2007strategies}. More recently, robotic on-orbit servicing has also been explored by DARPA in its Robotic Servicing of Geosynchronous Satellites (RSGS) program \cite{DARPARSGS} as well as by industries such as Northrop Grumman \cite{NorthropGrumann}. In addition, in-space manufacturing has also been explored by DARPA's Novel Orbital Moon Manufacturing, Materials, and Mass Efficient Design (NOM4D) Program \cite{DARPANOMAD} as well as other industry members such as Redwire.

In parallel to the hardware technologies for robotic ISAM, there is a significant need for logistics-driven approaches in the context of ISAM. The key questions that this paper will focus on include: 
\begin{itemize}
    \item A1Q1: how to analyze the performance of an ISAM architecture; 
    \item A1Q2: how to tactically plan and schedule the ISAM operations to satisfy the uncertain demands of the customer satellites under resource constraints; and
    \item A1Q3: how to strategically architect the ISAM infrastructure elements such as depots and vehicles.
    \end{itemize}
These questions can be answered through logistics-driven research as reviewed later.

\subsection{Application 2 (A2): Multi-Mission Space Exploration Campaigns}
Besides ISAM for satellites in Earth orbits, we are also interested in how space exploration campaigns can leverage in-space logistics infrastructure elements. The concept of an in-space propellant depot for space exploration has been studied extensively (e.g., \cite{Gaebler2009,Wilhite2012,Noevere2013,Ho2014-Acta}). In addition, the recent development of technologies for in-situ resource utilization (ISRU) on the Moon and Mars has become a game changer \cite{isru,shishko2017integrated,shishko2017mars,Kornuta2019}. For example, the successful oxygen generation on Mars by the Mars Oxygen In-Situ Resource Utilization Experiment (MOXIE) experiment has paved the way for Mars exploration with little reliance on Earth \cite{moxie}. Furthermore, reusable infrastructure elements and vehicles in space can also drastically change future space exploration. For example, NASA's Gateway concept is a representative example of in-space infrastructure for supporting exploration of the Moon and beyond \cite{Gateway}. To achieve sustainable deep space exploration, designing effective multi-mission space exploration campaigns leveraging such in-space logistics infrastructure elements is critical.

Like in the ISAM application, there is a significant need for logistics-driven approaches to designing and planning robotic and human space exploration campaigns. The key questions that this paper will focus on include: 
\begin{itemize}
    \item A2Q1: how to analyze the performance of a logistics strategy in the context of a space exploration campaign;
    \item A2Q2: how to plan and schedule the missions in a multi-mission exploration campaign with different types of trajectories (high-thrust, low-thrust, etc.) leveraging in-space infrastructure;
    \item A2Q3: how to design and size the exploration vehicles and resource infrastructure technologies such as depots and ISRU plants; 
    \item A2Q4: how to respond to uncertainties in launch delay, infrastructure performance, etc. in a space campaign; and 
    \item A2Q5: how to build the relationship between the governments and commercial players who supply in-space infrastructure.
\end{itemize} Logistics-driven research has played and will continue to play an important role in answering these questions.

\subsection{Application 3 (A3): Mega-Scale Satellite Constellations}
Another growing trend in low-Earth space is mega-scale satellite constellations. Many entities have proposed and developed such systems including OneWeb \cite{OneWeb}, SpaceX's Starlink \cite{SpaceX}, and Amazon's Project Kuiper \cite{Amazon}. The Department of Defense is also interested in concepts such as satellite swarms and proliferated low Earth orbit (pLEO). Although satellite constellations have been developed and launched in the past, including Iridium \cite{Iridium} and Globalstar \cite{Globalstar}, those conventional ones involved tens of satellites in each constellation. In contrast, those recent constellations will involve hundreds to thousands of satellites, which is a significantly larger scale than any conventional constellation system.

It is important to note that mega-scale constellations are not just a larger scale of small constellations. Instead, we encounter new logistics challenges due to the scalability issue of conventional methods. The key questions include: 
\begin{itemize}
    \item A3Q1: how to launch and deploy a mega-scale constellation when it does not fit in a single launch vehicle;
    \item A3Q2: how to analyze the system performance and allocate the on-orbit spares for the constellation so that it reliably satisfies performance requirements when failures would happen more frequently and in a more distributed way than conventional satellite systems;    
    \item A3Q3: how to flexibly reconfigure a satellite constellation when there is a change in the demand; and
    \item A3Q4: how to manage the commercial multi-stakeholder ecosystem for large-scale constellations.
\end{itemize} As the scale and complexity of the space systems of interest grow, logistics-driven techniques will be critical to address some of the inherent key challenges.

\subsection{Mapping of Literature to Each Application Question}
The following Table \ref{tab:mappingapplications} shows the mapping of each study reviewed in this paper to each key research question listed above. In our review, we primarily focus on the literature in recent years, but some earlier representative ones are also included to provide the context. The technical details of each work are revisited as they are re-categorized by the employed technical logistics-driven methods in Sec.\ref{methods}. Note that most of these questions are still open and new studies are being performed as we speak.
\begin{table}[h!]
    \centering
    \caption{Categorization of the State of the Art (SOTA) by Application Questions}
    \begin{tabular}{ll}
      \multicolumn{1}{c}{Application Questions} & \multicolumn{1}{c}{SOTA}\\
      \hline
      Application 1: ISAM for Satellites & \\
      \hline
     A1Q1: how to analyze the performance of an ISAM architecture & \cite{martin1988application,lamassoure2002space,saleh2002space,long2007orbit,Sears2018,Ho2020}\\
     A1Q2: how to tactically plan and schedule the ISAM operations & \cite{Bang2018,Verstraete2018,Bang2019,hudson2020,Jonchay2021,Jonchay2022,lee2023optimal}\\
     A1Q3: how to strategically architect the ISAM infrastructure elements such as depots and vehicles & \cite{Jonchay2021,Jonchay2022,shimane2024}\\
     \hline
      Application 2: Multi-Mission Space Exploration Campaigns & \\
      \hline
     A2Q1: how to analyze the performance of a logistics strategy & \cite{shull2006future,shull2006logistics,gralla2006modeling,de2007spacenet,siddiqi2009matrix,ishimatsu2010interplanetary,grogan2011matrix,grogan2011comparative,grogan2011space,moraguez2020benefits,capra2021spacenet}\\
     A2Q2: how to plan and schedule the missions in a multi-mission exploration campaign & \cite{Taylor2007,Arney2014,Ho2014-Acta,Ishimatsu2016,Ho2016,Chen2018,Jagannatha2018,Chen2019,Mcbrayer2019,Jagannatha2020}\\
     A2Q3: how to design and size the exploration vehicles and resource infrastructure technologies & \cite{Taylor2006,Chen2018,Chen2020,Chen2021-multi,Isaji2022}\\
     A2Q4: how to respond to uncertainties in launch delay, infrastructure performance, etc. & \cite{Chen2021,Takubo2022}\\
     A2Q5: how to build the relationship between the governments and commercial players & \cite{grogan2012federated,shishko2018affordable,sarton2020commercial,Chen2022}\\
     \hline 
     Application 3: Mega-Scale Satellite Constellations & \\
      \hline
     A3Q1: how to launch and deploy a mega-scale constellation & \cite{de2004staged,Lee2018,sung2023optimal}\\
     A3Q2: how to analyze the system performance and allocate the on-orbit spares for the constellation & \cite{dishon1966communications,kelley2004minimizing,Zheng2009,Jakob2019}\\
     A3Q3: how to flexibly reconfigure a satellite constellation when there is a change in the demand & \cite{de2008optimal,Lee2023regional}\\
     A3Q4: how to manage the commercial multi-stakeholder ecosystem for large-scale constellations & \cite{grogan2014multi,grogan2015interactive,grogan2016multi,Guzzetti2022satellite,Qureshi2022modeling,Qureshi2023table-top}
    \end{tabular}
    \label{tab:mappingapplications}
\end{table}
\newpage
\section{Categorization of State of the Art by Logistics-Driven Methods}
\label{methods}
This section categorizes the same set of literature reviewed above by the logistics-driven methods used in each study to aid the researchers in understanding the state of the art from the technical perspective.

Terrestrial logistics has a long history of research and practice, from which we have learned a lot to advance space logistics research. This subsection summarizes three key method categories, each of which covers a subfield of logistics research that is particularly relevant to space applications: (1) Network Flow Modeling and Optimization for Logistics Planning and Scheduling; (2) Probabilistic Modeling and Queueing Theory for Logistics Performance Analysis; and (3)  Inventory Control for Resource Infrastructure Operations Management. 
Each subfield has been largely developed with terrestrial applications in mind, and thus the developed theories and techniques cannot be directly applicable to space applications. The challenges of applying terrestrial logistics-driven techniques are discussed for each subfield, as well as how the literature has addressed some of these challenges.
Lastly, the set of literature that extends beyond these three key method categories in response to space-unique challenges is discussed.
For more technical details on each terrestrial logistics subfield, refer to the textbooks in the Operations Research (OR) field like Ref. \cite{LarsonOdoni,SimchiLeviChenBramel}.


\subsection{Method 1 (M1): Network Flow Modeling and Optimization for Logistics Planning and Scheduling}
\label{subfield1}
\subsubsection{Motivation and Overview}
\begin{figure}[h]
\begin{center}
\includegraphics[width=\textwidth]{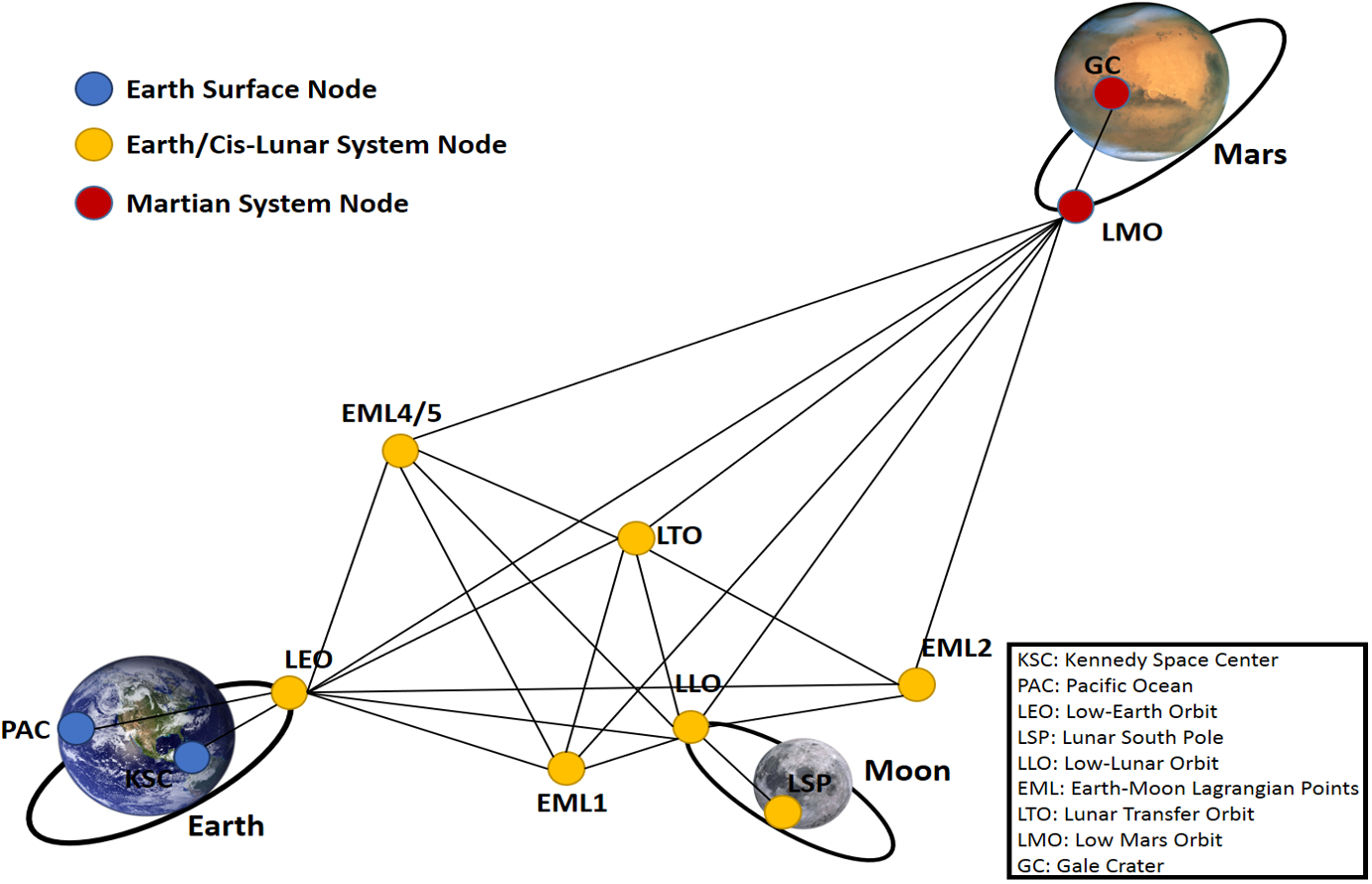}
\caption{Example of network modeling for Earth-Moon-Mars exploration logistics \cite{Chen2019}.}
\label{fig1}
\end{center}
\end{figure}

As we consider how to manage a system of spacecraft and infrastructure distributed over the vast space in Earth's orbits, cislunar, and interplanetary space, it is natural to draw an analogy with how terrestrial logistics operations manage their fleet of vehicles and commodity flows, for which network modeling and optimization has been a primary tool. Figure \ref{fig1} shows an example of network modeling for Earth-Moon-Mars exploration logistics.

A network-based approach models the system as a set of nodes and the arcs connecting them and employs various techniques to analyze and optimize the commodity flow over that network. For example, for a space exploration campaign problem, each orbital staging point (e.g., geosynchronous orbits, Lagrangian points, low-lunar orbits) is modeled as a node and the trajectories connecting them are modeled as arcs. For an on-orbit servicing problem, the nodes can also correspond to the customer satellites to be serviced. The commodities include everything that flows over the network, including the propellant, spares, tools, payload, and vehicles themselves. The typical goal for such a problem is to find a vehicle/commodity transportation routing plan and schedule to minimize a cost function, which can be defined based on total campaign launch mass, cost, or other relevant metrics. Some of the well-known problem types in network optimization include the traveling salesman problem (TSP), which finds the optimal route to visit all nodes once and return to the origin node, and its generalized variant, the vehicle routing problem (VRP), which finds the optimal set of routes for a fleet of vehicles to deliver to a given set of customers. 
Another relevant class of network problem is the facility location problem (FLP), which finds the optimal location of the facilities such as depots and servicing stations to optimally serve the distributed customers. 
Thus, we can draw analogies between a spacecraft mission design problem and a TSP/VRP, and between the propellant depot or ISAM facility location problem with an FLP.

\subsubsection{Challenges in Space Logistics Applications}
Unfortunately, optimizing and analyzing space logistics missions is substantially more complex than just solving a conventional TSP/VRP/FLP
 over a $\Delta$V map of orbits. 
There are some unique challenges that prevent a direct application of existing terrestrial logistics approaches.

The most obvious difference between terrestrial logistics and space logistics is the orbital mechanics. Trajectory beyond Earth's orbits cannot be modeled as a two-body problem and thus, the generation of the $\Delta$V or TOF will require a large computational effort. For high-thrust trajectories, these quantities can be pre-computed (e.g., via Lambert problem solvers) before optimizing the commodity flow, but this approach is not feasible for low-thrust trajectories, where the $\Delta$V and TOF can be flow-dependent; namely, if we assume a constant-thrust engine, the acceleration (and thus the $\Delta$V and TOF) depends on the mass of the commodities that the spacecraft carries. Thus, the logistics optimization problem and the trajectory design problem are coupled. 
In addition, particularly when various trajectory options are considered for logistics missions in cislunar space and beyond, drastically different time scales of the flight (from hours to months or years) would need to be included in the same network, making the modeling of the time dimension challenging. On top of that, the orbital mechanics also lead to complex time-dependent trajectory performance; for example, the relative locations of the depots and servicers could be time-dependent, and thus cannot be pre-computed as a ``map.'' These challenges are unique to space applications and make the conventional network optimization formulation ineffective or infeasible.

Besides orbital mechanics, the nature of space infrastructure also poses additional unique challenges. In terrestrial applications, when an infrastructure/facility is built and used, the deployment phase is typically relatively short compared to the utilization phase of the infrastructure. However, that is not the case for space. As can be seen in examples like the International Space Station (ISS), the deployment and assembly phase of infrastructure in space can be nearly as long as its utilization phase; the boundary between deployment and operation is blurred. With ISRU or in-space manufacturing, we can even use the resources generated by the previous stage of the infrastructure to deploy more infrastructure (this concept is sometimes referred to as bootstrapping strategy \cite{Kornuta2019}). Thus, we need to consider the deployment of infrastructure together with its operations. 

\subsubsection{State of the Art}
\label{subfield1review}
Fortunately, many (yet not all) of the above challenges have been tackled over 30 years of research in space logistics.
Network modeling and optimization are the most well-explored fields in space logistics modeling; see Fig. \ref{fig1} as an example. One of the first projects that leveraged network modeling for space mission design was SpaceNet \cite{de2007spacenet}, which is a simulation and visualization software for space logistics operations. This project and its subsequent studies (e.g., Refs. \cite{gralla2006modeling,grogan2011comparative,grogan2011space,siddiqi2009matrix,grogan2011matrix,ishimatsu2010interplanetary,capra2021spacenet}\footnote{Some of the cited works (e.g., Ref. \cite{siddiqi2009matrix}) did not necessarily leverage the network structure explicitly, but they are included as they leveraged the network modeling framework in SpaceNet.}) primarily focused on supporting human space exploration campaigns with the then-active Constellation Program in mind. Since then, network and graph-theoretic modeling have played a key role in analyzing complex space logistics missions. Some of the representative works that focused on network modeling for space logistics missions include Refs. \cite{Arney2014,Ishimatsu2016,Ho2014-Acta,Ho2016}. A more recent work specifically focused on the depot location problems for ISAM applications with a facility location problem formulation \cite{shimane2024}. A variant of network flow models, the assignment problem, has also been used for satellite constellation deployment and reconfiguration analysis \cite{de2008optimal,sung2023optimal,Lee2023regional}.

From the optimization perspective, a variety of different formulations have been proposed to efficiently and effectively optimize the space exploration and/or ISAM mission design process. While some studies used metaheuristics such as genetic algorithms or particle swarm optimization \cite{Verstraete2018,Jagannatha2018} or simply evaluating all permutations of missions \cite{hudson2020}, most studies developed network-based mixed-integer linear programming (MILP) formulations so that commercial solvers such as Gurobi can solve them within guaranteed optimality gap. Note that, although the trajectory problem is a nonlinear problem, a two-phase approach has been proposed to decouple the nonlinear trajectory problem from the integer linear network problem \cite{Bang2018,Bang2019,lee2023optimal,sung2023optimal}. Two major types of network formulations for logistics optimization include (1) a path formulation and (2) a node-arc formulation. A path formulation defines its decision variables as the commodity flows over each possible vehicle's and commodity's path over the network, whereas a node-arc formulation defines its decision variables as the vehicle/commodity flows over each arc. The choice of formulation depends on the problem's structure, constraints, and the solution method's characteristics; for example, while a path formulation would be effective in modeling a generic path-based cost and/or feasibility conditions, a node-arc formulation would be more effective when the vehicle's or commodity's path is not easily definable or enumerable (e.g., reconfigurable vehicles; ISRU resources).
In the literature, path formulations have been typically developed in conjunction with column-generation techniques; such examples include Martian surface vehicle routing \cite{Ahn2012} and lunar mission vehicle routing (with no ISRU/ISAM) \cite{Taylor2007}. On the other hand, arc-based formulations have been developed when in-space infrastructure elements are involved, such as multi-mission Mars-Moon-asteroid space campaign design with ISRU \cite{Ishimatsu2016,Ho2014-Acta,Ho2016}. Recently, another formulation, a path-arc formulation, was developed, which leverages the path formulation for the vehicle and the node-arc formulation for the commodities \cite{Mcbrayer2019}.

To tackle the trajectory-related challenges, attempts have been made to develop approximate surrogate models to compute the low-thrust, low-energy trajectory effectively leveraging Q-law and 3-body dynamics \cite{Jagannatha2020}. Such approximation models can be incorporated into MILP optimization through piece-wise linear approximation so that both (approximate) trajectory design and logistics optimization are performed concurrently.

To address the time-related challenges, the literature has also developed multiple techniques to model and optimize the dynamic network flow. A classical method to consider the dynamic network flow is to use a time-expanded network, where the nodes are copied over the time dimension with predefined time steps and network flow is optimized over that expanded network. However, due to the significant time-scale differences in the problem and sometimes unknown TOF without solving the logistics problem (i.e., in the low-thrust trajectories), a classical time-expanded network can be ineffective or infeasible. To tackle this challenge, multiple variants of time-expanded networks have been developed. For the case where the different time scales are known beforehand, a bi-scale time-expanded network was developed \cite{Ho2014-Acta}. For a specific type of problem such as Mars exploration leveraging the lunar resources, a partially-static time-expanded network was also proposed as an approximation technique \cite{Ho2014-Acta}. For a campaign that will largely repeat once entering a steady state (e.g., resupply logistics for Mars habitats), a partially-periodic time-expanded network was developed to effectively concurrently optimize both the build-up of the space system and its steady operations \cite{Chen2019}. For the flow-dependent TOFs, where any of the above variants of time-expanded networks would not be applicable, an event-driven network was developed to optimize the logistics flow without specifying the time steps \cite{Jagannatha2020}.

To consider the unique resource transformation mechanisms with ISRU and ISAM infrastructures, a new network formulation was developed in the space logistics community: generalized multi-commodity network flow (GMCNF), including its static \cite{Ishimatsu2016} and dynamic versions \cite{Ho2014-Acta}. This formulation generalizes the conventional multi-commodity flow formulation so that it can model the commodity transformation along the arcs (including resource generation or consumption); this general capability has enabled the modeling and optimization of missions with a significantly larger variety of space infrastructure systems such as ISRU \cite{Ishimatsu2016,Ho2014-Acta} and ISAM systems \cite{Jonchay2021,Jonchay2022}. Combined with the partially periodic time-expanded network mentioned above \cite{Chen2019}, we can also model the deployment of resource infrastructure elements and their operations utilization concurrently.

\subsection{Method 2 (M2): Probabilistic Modeling and Queueing Theory for Logistics Performance Analysis}
\label{subfield2}
\subsubsection{Motivation and Overview}
\begin{figure}[h]
\begin{center}
\includegraphics[width=\textwidth]{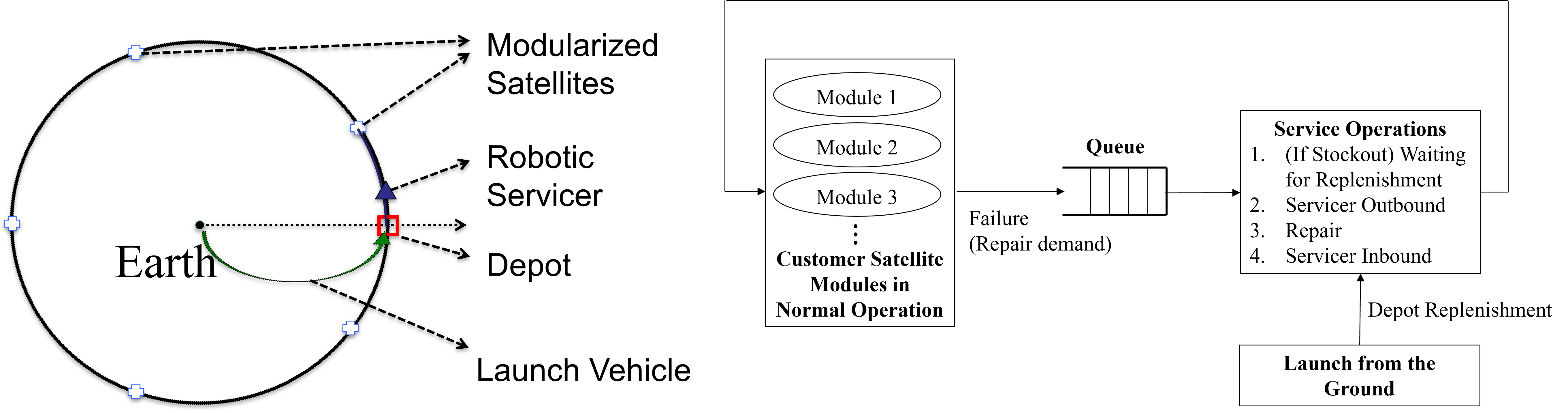}
\caption{Example of queueing modeling for an on-orbit servicing architecture \cite{Ho2020}.}
\label{fig2}
\end{center}
\end{figure}

As we consider how to analyze the logistics performance of a space system under uncertainties, like a multi-satellite system under geometrically distributed uncertain failures and a servicing infrastructure to respond to them, we can learn from how fire stations or ambulance systems are designed to respond to uncertain demands distributed in a city. The probabilistic modeling and spatial queueing theory have been very useful tools to that end. Figure \ref{fig2} shows an example of queueing modeling for an on-orbit servicing architecture.

Queueing theory is a mathematical theory to analyze the performance of a queue. With the service time distribution and demand rate information as inputs, queueing theory analyzes the waiting time of a randomly arriving customer as well as how long the queue is expected to be. The simplest queueing model, M/M/1, assumes a Poisson demand arrival rate and an exponentially distributed service time and uses a Markov chain to model the system; more advanced models have been developed such as a queue with multiple servers (e.g., M/M/m), a queue with capacity (e.g., M/M/1/c), 
a queue with general service time distribution (e.g., M/G/1), and a queue with multiple classes of customers with different priorities (priority queues). To model fire station or ambulance systems, spatial queueing models have been used, which are based on the M/G/1 model with the service time distribution built based on the geometric relationship between the servicer(s) and the potential customers \cite{LarsonOdoni}. Thus, analogously, to analyze the service responsiveness of an in-space servicing system, for example, we can consider a similar spatial queueing model (i.e., \textit{orbital} queueing model) with the service time distribution derived based on the orbital relationship between the servicer(s) and the potential customers.

\subsubsection{Challenges in Space Logistics Applications}
Developing an \textit{orbital} queueing model will involve challenges beyond a direct application of M/G/1 models. Each of these challenges is not necessarily unexplored in terrestrial applications, but the combination of them, coupled with orbital mechanics, makes the space application challenging.

One difference between the ISAM facility servicing and fire station or ambulance systems is the limited resources. For any refueling and servicing operations, we need fuel and spares launched from the ground. Thus, we need to consider an inventory management strategy for the storage space of the servicer and the depot in orbit if there is one. (The study of inventory management is discussed in Sec. \ref{subfield3}.) This coupling between the queueing model and inventory management makes the analysis challenging.

Additionally, the number of satellites a servicer serves is relatively limited. Thus, the demand distribution of the satellite population would depend on the number of satellites that have failed and are being repaired (i.e., the less the number of functioning satellites, the lower the failure rate for the entire system.); note that this feature can often be ignored for terrestrial applications because of the large population that a fire station or an ambulance system serves. A queue with a finite number of customers (often referred to as a finite-source queue) has been studied \cite{sztrik2001finite}, but its analysis is a substantially more challenging problem than conventional queues, particularly with generalized assumptions about the service time, the number of servicers, the classes of customers with priorities, etc.

Furthermore, the dynamic nature of the demand and service rates can make the problem challenging. For example, the demand rate (e.g., satellite failure rate) can be different depending on the phase of the satellite's lifetime. In addition, the required service time and cost are also time-dependent if the relative position of the servicing depots and the customer satellites changes over time. An effective model needs to address these challenges so the system can be useful for realistic space logistics analysis. 

\subsubsection{State of the Art}
\label{subfield2review}
Modeling in-space logistics using probabilistic modeling techniques and queueing theory has been studied since the 1980s \cite{martin1988application}, but their theoretical developments for realistic space applications have been relatively limited.

Earlier studies focused on modeling ISAM applications based on queueing theory and analyzed them with discrete-event simulations or agent-based simulations. Simulations have been used for on-orbit servicing analysis from the early 2000s \cite{long2007orbit}. More recently, attempts were made to model on-orbit servicing, recycling, and manufacturing applications using discrete-event simulations \cite{Sears2018,moraguez2020benefits}, in which the impacts of having different on-orbit recycling and manufacturing technologies on the performance of on-orbit operations were quantitatively evaluated.

While simulations can tackle general scenarios, they are stochastic and computationally expensive, and thus not effective for analyzing and optimizing the system; rather, there is a high need for analytic modeling of space logistics systems. Various studies have employed probabilistic modeling techniques based on a Markov model for analyzing the system availability of a constellation and designing its spare strategies \cite{Zheng2009,kelley2004minimizing} or its servicing strategies \cite{saleh2002space,lamassoure2002space}. A more recent study developed a semi-analytic model for an on-orbit servicing system based on the finite-source M/G/1 queueing model (i.e., M/G/1/K/K model) \cite{Ho2020}. In this study, to model the limited launch resupply opportunities for the spares from the ground, an inventory control technique was also integrated (see Sec.\ref{subfield3review}). The resulting analytic model matches well with the simulation model and can cut down the computational time from hours/days for simulations down to seconds, enabling large-scale tradespace exploration and optimization in the early stage of the design.

\subsection{Method 3 (M3): Inventory Control for Resource Infrastructure Operations Management}
\label{subfield3}
\subsubsection{Motivation and Overview}
\begin{figure}[h]
\begin{center}
\includegraphics[width=\textwidth]{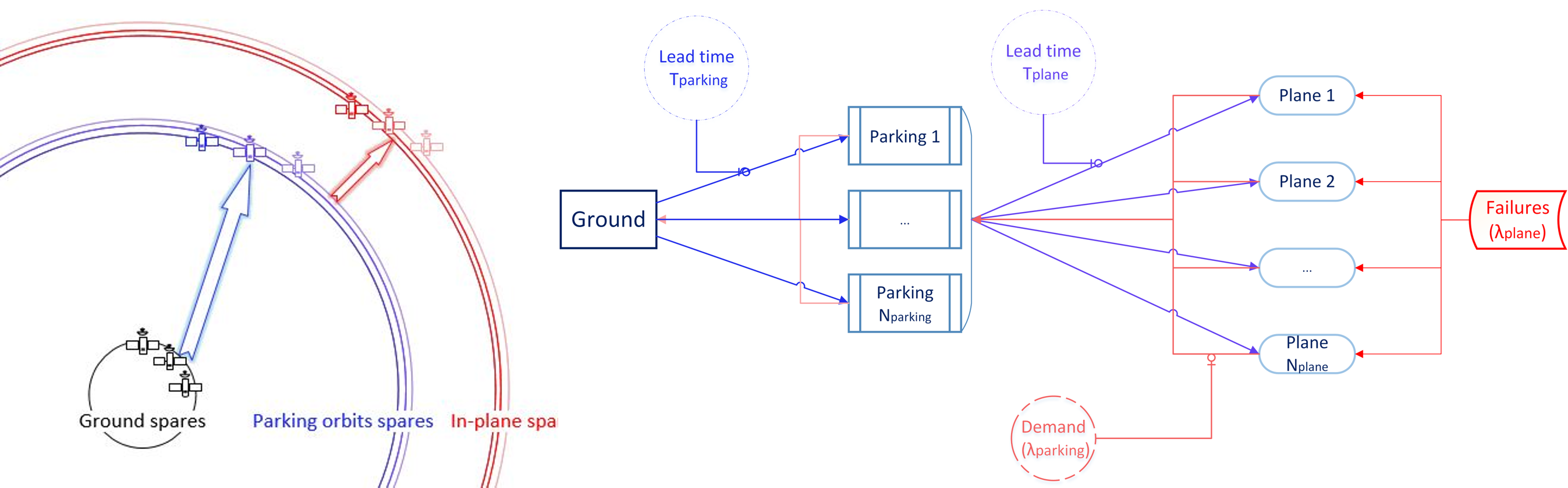}
\caption{Example of inventory modeling for spare strategy analysis for a mega-scale satellite constellation \cite{Jakob2019}.}
\label{fig3}
\end{center}
\end{figure}

As ISAM becomes a trend in space development, it is critical to explore the concepts of in-space propellant depots, spare/tool warehouses, and other types of infrastructure. These infrastructure elements need to be regularly refueled or resupplied from the ground (or possibly ISRU facilities on the Moon or Mars) and thus, the operations management of the inventory can become a research question. To this end, we can draw an analogy between the inventory control of in-space depots/warehouses and that of terrestrial depots/warehouses. Figure \ref{fig3} shows an example of inventory modeling for spare strategy analysis for a mega-scale satellite constellation.

Inventory management is a subfield in logistics that analyzes and optimizes inventory policies, including the timing and quantity of orders, to minimize cost. Typically, the cost metric is the summation of the holding cost, purchase cost, fixed ordering cost, and shortage penalty. Inventory management becomes different when there are uncertainties in the demand and the lead time (i.e., the time between the order and its delivery). Given a cost metric, we can formulate an optimization problem to find the optimal policy under uncertainties; one such problem is the Newsvendor model \cite{SimchiLeviChenBramel}. When there are constraints to the system on how the orders need to be made (e.g., periodic launch, launch capacity), we can develop a parametric inventory policy and optimize their parameters. Some popular policies include (1) a periodic review, order-up-to policy, which periodically orders to fill the inventory up to a certain level; (2) a reorder point, order quantity policy, which orders a fixed quantity as soon as the inventory level drops below a threshold; and (3) a reorder point, order-up-to policy, which orders to fill the inventory up to a certain level  as soon as the inventory level drops below a threshold. These parametric policies can be applied and optimized to space logistics context as well depending on the structure of the problem. Furthermore, if we can consider multiple layers of inventory (e.g., manufacturer $\rightarrow$ warehouses $\rightarrow$ retailers), we can consider a multi-echelon inventory model; this can be useful to spacecraft spare inventory, where Earth is modeled as manufacturer, spare orbits/depots as warehouses, and customer satellites as retailers. 

\subsubsection{Challenges in Space Logistics Applications}
Applying inventory management to space logistics applications involves orbital mechanics, which often causes additional challenges in the modeling of stochastic demand and lead time. 

For the satellite refueling or repairing applications in the context of ISAM, the demand corresponds to the customer satellites' needs for fuels, materials, and spares. As reviewed in Sec. \ref{subfield2}, the demand depends on the phase of the satellites' lifetime and thus can be time-varying. Furthermore, if servicing operations are involved, the demand is coupled with the transportation of the servicer (governed by the spatial queueing model). Thus, an integrated model between queueing and inventory management is often needed.

In addition, a more unique challenge in space logistics application is the lead time distribution. The lead time corresponds to the time between the order and its delivery, which is modeled based on orbital mechanics. However, due to the dynamic nature of orbital mechanics, the relative position of the depots and customer satellites can vary over time. Thus, when multiple depots and multiple customer satellites are considered, the satellites need to be resupplied from the closest depot with a sufficient inventory level. This additional complexity in the modeling, which can be challenging if complex orbital mechanics is involved, is one of the unique requirements for space applications.

\subsubsection{State of the Art}
\label{subfield3review}
The idea of inventory management has been conceived for satellite replenishment since the 1960s \cite{dishon1966communications}. When space logistics research was applied to human space exploration, inventory management was also recognized as an important area \cite{shull2006logistics,shull2006future}. However, it was not until recently that mathematical inventory control theory was rigorously applied to tackle realistic large-scale space logistics problems. 

One of the recent studies that focused on developing an inventory control model is Ref. \cite{Jakob2019}. In this work, a multi-echelon inventory model was developed for spare strategies for mega-scale constellations. In this work, separate spare parking orbits are considered at a lower altitude than the original constellation orbit, and an analytic multi-echelon inventory model based on a reorder point, order quantity policy was developed by modeling Earth as the manufacturer, spare parking orbits as warehouses, and the constellation as the retailers. The resulting inventory model can be optimized efficiently to identify the optimal number of parking orbits and their locations as well as their inventory control strategies.

Another work that optimized the inventory policy under uncertainties is Ref. \cite{Chen2021}. In this work, a model was developed to optimize space station logistics under launch schedule uncertainties. Specifically, the inventory policy (i.e., referred to as a decision rule) was parameterized by the safety stock level of the onboard inventory, and this parameter was optimized along with the logistics decisions to balance the loss due to supply shortage and the cost of extra supply (i.e., safety stock). Namely, it optimized the operational strategy to respond to the uncertainties in the launch delays.

In addition, as reviewed earlier, inventory control was also used in combination with queueing models in Ref. \cite{Ho2020}. In this work, the inventory model is used to estimate how often a stockout happens (i.e., not enough spares are available for further servicing missions) and how much delay is expected when it happens. This stochastic model is used as part of the service time distribution model, which is then fed back into the queueing model to evaluate the performance of the ISAM system.

\subsection{Extensions of Basic Models}
Beyond the above basic models introduced, the literature has explored the extensions of these basic methods to address unique challenges in space logistics challenges. Some are reviewed below.

\subsubsection{Concurrent Design of Vehicle/Infrastructure and Space Logistics Network}
One unique feature of space missions is that we often design space vehicles and infrastructure elements dedicated to a certain space campaign. Thus, we also have an interesting research question in space logistics regarding how to achieve the concurrent design optimization of the vehicle, infrastructure, and logistics commodity flow network. This problem often involves nonlinearity and multi-disciplinary design analysis due to the nature of the spacecraft and space infrastructure design. Thus, to solve such a problem, we need new strategic formulations and solution methods beyond traditional logistics research.

To enable the incorporation of vehicle sizing and infrastructure design, a few techniques have been developed depending on the required fidelity. If a linear vehicle and infrastructure sizing model is sufficient, a regular MILP formulation can be used \cite{Ho2014-Acta,Ho2016,Ishimatsu2016}. A higher-fidelity version would be to use a more piece-wise linear approximation for the vehicle sizing, which can still be converted into a MILP formulation but requires more variables and constraints \cite{Chen2018}. For a yet higher-fidelity case where a nonlinear model needs to be used, an embedded optimization \cite{Taylor2006} was developed earlier, and a more advanced and effective Augmented Lagrangian Coordination (ALC) approach has been developed recently \cite{Isaji2022}. The ALC-based approach is based on multi-disciplinary optimization and leverages the unique structure of the problem by decomposing the whole mixed-integer nonlinear programming problem into a set of mixed-integer quadratic programming problems and nonlinear programming problems, each of which can be solved using specialized solvers. Finally, to incorporate the multi-subsystem interaction within the ISRU or other in-space infrastructure into the logistics optimization, a multi-fidelity network flow was developed to capture these interactions effectively \cite{Chen2020,Chen2021-multi}.

\subsubsection{Incorporating Operational Uncertainties in Space Logistics Network Optimization}
Although earlier space logistics research has primarily focused on deterministic problem setting, recent studies started to incorporate uncertainties such as launch delays, demand fluctuation, and infrastructure performance decay. 

A major approach is to characterize the uncertainties into a set of possible scenarios and optimize the strategy to respond to these scenarios together with the baseline design. The most straightforward formulation to this end is stochastic optimization, which has been applied to multi-stage mega-scale constellation deployment under demand uncertainties \cite{de2004staged,Lee2018}. This approach is simple but cannot effectively provide a strategy to respond to the uncertainties. As an alternative, as reviewed earlier, the decision-rule-based optimization was developed where the strategies (i.e., decision rules) to respond to the uncertainties are parametrically defined and optimized (e.g., safety stock to prepare for launch delays) \cite{Chen2021}. Furthermore, to handle more generic uncertainties in a space campaign, where the decision rules cannot be intuitively parameterized, a Hierarchical Reinforcement Learning (HRL) framework was developed \cite{Takubo2022}. Leveraging the strength of reinforcement learning (RL) to handle uncertainties and MILP's strength to handle a large variable space well, the HRL framework uses RL to provide high-level guidance and uses MILP to optimize the detailed action to achieve that guidance. This success was one of the first studies that incorporated machine learning into space logistics. 

Another approach to determining how to respond to uncertainties is to reoptimize as soon as we learn more information about the uncertainties \cite{Verstraete2018,hudson2020}. Building upon this idea, a rolling horizon approach was developed for ISAM applications to respond to newly obtained demand information over time \cite{Jonchay2021,Jonchay2022} and a reconfiguration approach was developed for satellite constellation configuration modification as we acquire new demand \cite{de2008optimal,Lee2023regional}.

\subsubsection{Space Logistics Modeling for Commercialization}
Most existing works have focused on the centralized optimization of the entire space logistics network. However, given the recent trend of space commercialization, this is not necessarily the best modeling approach. Thus, an interesting recent extension of the above work is to incorporate multiple players in space logistics. To this end, Ref. \cite{Chen2022} extended network optimization research with the game theory to optimize the incentive that the government should provide to the commercial players (e.g., infrastructure providers). Other works used alternative approaches to examine the affordability and commercialization of space missions \cite{shishko2018affordable} and investigated a method to use the commercial suitability as an explicit criterion of the initial architecture selection \cite{sarton2020commercial}. Furthermore, various studies developed an interactive game as well as a simulation framework to assess the multi-stakeholder space systems \cite{grogan2012federated,grogan2014multi,grogan2015interactive,grogan2016multi,Guzzetti2022satellite,Qureshi2022modeling,Qureshi2023table-top}.

\section{Summary of Reviewed Literature}
\label{sec:summary}
Table \ref{tab:mapping} summarizes how each reviewed literature is mapped to the corresponding application area and logistics-driven method. Note that some papers span over multiple methods (e.g., Ref. \cite{Ho2020}). While not all combinations of the application areas and logistics-driven methods have an equal amount of room for exploration, this table shows some well-explored and under-explored combinations of application areas and methods.

\begin{table}[h!]
    \centering
    \caption{Mapping of the Applications and Methods}
    \begin{tabular}{lllll}
      Applications/Methods & M1 & M2 & M3 & Extensions  \\
      \hline
      A1: ISAM for Satellites & \cite{Jonchay2021,Verstraete2018,hudson2020,Jonchay2022,Bang2018,Bang2019,lee2023optimal,shimane2024} & \cite{martin1988application,long2007orbit,lamassoure2002space,saleh2002space,Sears2018,Ho2020} & \cite{Ho2020} & \cite{Jonchay2021,Jonchay2022,Verstraete2018,hudson2020}\\
      A2: Space Exp. Campaigns & \cite{Ho2014-Acta,gralla2006modeling,de2007spacenet,siddiqi2009matrix,grogan2011matrix,ishimatsu2010interplanetary,grogan2011comparative,grogan2011space,capra2021spacenet,Taylor2007,Arney2014,Ishimatsu2016,Chen2019,Jagannatha2018,Jagannatha2020,Mcbrayer2019,Ho2016} & \cite{moraguez2020benefits} & \cite{shull2006future,shull2006logistics,Chen2021} & \cite{Ho2014-Acta,Ho2016,Ishimatsu2016,Chen2018,Chen2020,Chen2021-multi,Taylor2006,Isaji2022,Chen2021,Takubo2022,Chen2022,grogan2012federated,sarton2020commercial,shishko2018affordable}\\
     A3: Satellite Constellations & \cite{de2008optimal,sung2023optimal,Lee2023regional} & \cite{Zheng2009,kelley2004minimizing} & \cite{dishon1966communications,Jakob2019} & \cite{de2004staged,Lee2018,de2008optimal,Lee2023regional,grogan2014multi,grogan2015interactive,grogan2016multi,Guzzetti2022satellite,Qureshi2022modeling,Qureshi2023table-top}\\
    \end{tabular}
    \label{tab:mapping}
\end{table}

\begin{figure}[h!]
    \centering
    \begin{subfigure}{\linewidth}
        \input{tikz_figs/before_2009}
    \end{subfigure}
    \begin{subfigure}{\linewidth}
        \input{tikz_figs/2010s}
    \end{subfigure}
    \begin{subfigure}{\linewidth}
        \input{tikz_figs/2020s}
        \vspace{-10pt}
    \end{subfigure}
    \caption{Number of citations (as of Feb. 20, 2024; according to Google Scholar) of the reviewed papers in each decade. }
    \label{fig:mapping_timeline}
\end{figure}
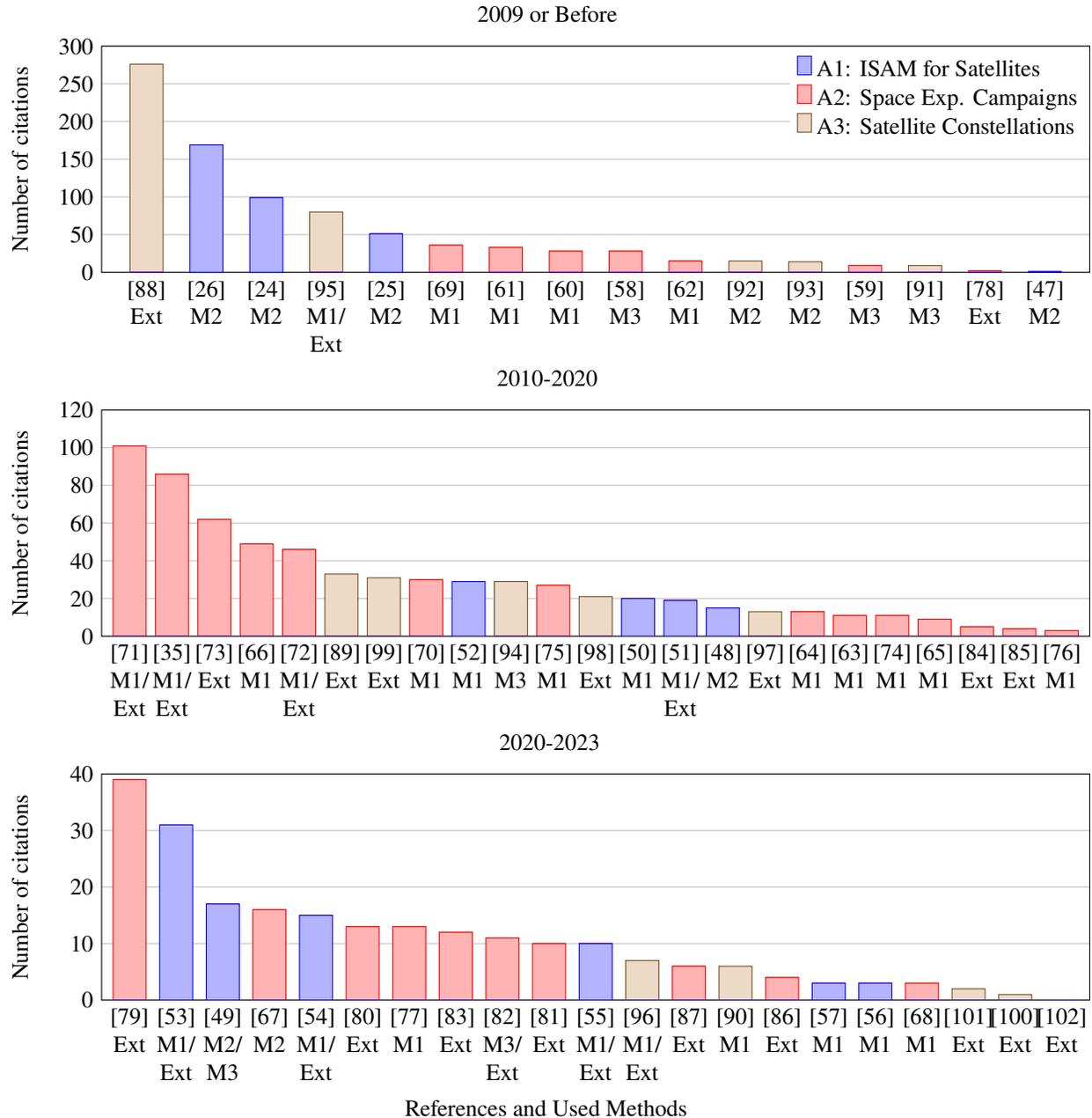

In addition, Figure \ref{fig:mapping_timeline} shows the number of citations (as of Feb. 20, 2024; according to Google Scholar) for each reviewed paper. The data is split into three time frames, the 2000s or earlier, the 2010s, and the 2020s, to show the most impactful papers in each decade in the space logistics area. Each reference also has its application (A1/A2/A3) and method(s) (M1/M2/M3/Ext) listed, where Ext stands for one of the extensions. Although the number of citations does not necessarily indicate the value of the work, this figure shows the trend and the shift of focus over time in terms of both the applications and methods.

\captionsetup[subfigure]{skip=-10pt}

While many questions have been addressed by space logistics research, there is enormous room for growth. Some examples of the open fundamental research questions are: (1) scalable algorithms (including heuristics) for large-scale space logistics network optimization: (2) complex queueing networks for space logistics analysis; and (3) dynamic multi-item inventory for space logistics applications, among others. Furthermore, given the nature of the growing space logistics applications, we can expect to see logistics-driven methods applied to more problems beyond the ones reviewed above. Drawing analogies between more terrestrial logistics applications and space logistics applications can lead to new research questions. Finally, most space logistics works have assumed simplified trajectory, vehicle, and operational models to enable large-scale tradespace exploration and optimization; once the solutions are obtained from these low-/mid-fidelity models, we need higher-fidelity simulations to validate these results. Such validations are critical yet largely unexplored research areas.

\section{Conclusion}
\label{sec:conclusion}
This paper provides an overview of the state of the art for space logistics modeling and optimization. The recent literature is categorized in two ways: (1) by application questions that are addressed; and (2) by logistics-driven methods that are used in the studies. The applications reviewed in this paper are determined by the needs of the community; the specific applications of interest include: (1) ISAM for satellites; (2) multi-mission space exploration campaigns; and (3) mega-scale satellite constellations. The logistics-driven methods reviewed in this paper are determined by mapping the relevant major logistics research subfield to space applications; the connection between each logistics research subfield to space logistics applications is reviewed and the relevant unique challenges in space applications are discussed. The reviewed logistics-driven methods include: (1) network flow modeling and optimization for logistics planning and scheduling; (2) probabilistic modeling and queueing theory for logistics performance analysis; and (3) inventory control for resource infrastructure operations management. Several extensions of these basic models that tackled unique challenges in space applications are also reviewed. We expect that the first application-based categorization helps practitioners understand the existing research that can answer their application questions, whereas the second method-based categorization helps the researchers to understand the technical perspective of the literature and identify new under-explored research directions.

\section*{Acknowledgment}
The author thanks the members of the AIAA Space Logistics Technical Committee for their helpful comments. The author also thanks Masafumi Isaji and Yuri Shimane for their suggestions about data visualization and improvement of the paper. While the initial version of the work was not associated with a funding source, its refinement and revision was conducted with support from the Air Force Office of Scientific Research
(AFOSR), as part of the Space University Research Initiative (SURI), under award number FA9550-23-1-0723.

\section*{Author Biography}
\begin{figure}[h]
\begin{center}
\includegraphics[width=0.3\textwidth]{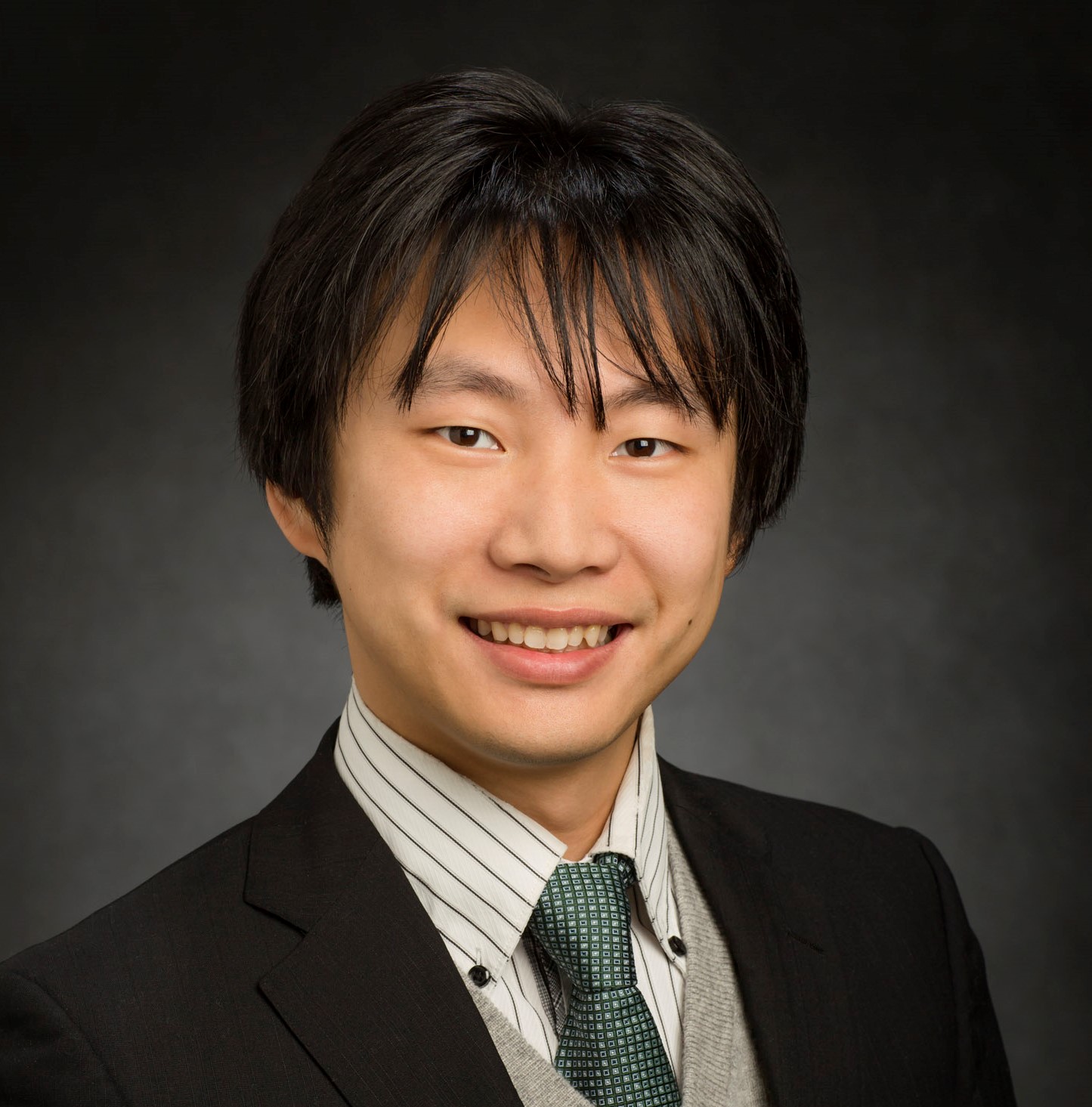}
\end{center}
\end{figure}
Dr. Koki Ho is the Dutton-Ducoffe Professor, Associate Professor, and the director of the Space Systems Optimization Group in the Daniel Guggenheim School of Aerospace Engineering at Georgia Institute of Technology. His research focuses on developing modeling and optimization methods for rigorous space mission analysis and design. His unique research connecting logistics-based modeling, optimization, systems engineering, and space applications has provided a substantial impact on modern and future space missions that involve multiple missions, multiple vehicles, and reusable infrastructure elements. Dr. Ho earned his Ph.D. at the Massachusetts Institute of Technology and his bachelor’s and master’s degrees at the University of Tokyo. He is the recipient of the NSF CAREER Award (2020), the NASA Early Career Faculty Award (2019), the DARPA Young Faculty Award (2019), and the Luigi Napolitano Award (2015), and he is a co-author of one of the most downloaded Acta Astronautica articles. Dr. Ho is the Chair of the AIAA Space Logistics Technical Committee (2017–2024) and also serves on the Steering Committee for the NASA-funded Consortium for Space Mobility and ISAM Capabilities (COSMIC). 

\bibliography{sample}

\end{document}

%% file: tikz_figs/before_2009.tex

\def\filename{data.csv}

\pgfplotstableread[col sep=comma]{data.csv}\datatable
\def\figtitle{2009 or Before}
\begin{tikzpicture}
\begin{axis}[
    ybar stacked,
    bar width=.5cm,
    width=\textwidth,
    height=5cm,
    xtick=data,
    xticklabels from table={\datatable}{x_label},
    x tick label style={
        align=center, 
        },
    xtick style={draw=none}, 
    enlarge x limits=0.05, 
    ylabel={Number of citations},
    ymin=0,
    ymax=300,
    ytick={0,50,100,150,200,250,300},
    ymajorgrids, 
    ytick style={draw=none}, 
    legend style={
        at={(1,0.98)}, 
        anchor=north east, 
        legend columns=1, 
        draw=none, 
        },
    legend cell align={left},
    ]

\addplot+[] table [x expr=\coordindex, y=A1, col sep=comma] {\datatable};
\addplot+[] table [x expr=\coordindex, y=A2, col sep=comma] {\datatable};
\addplot+[] table [x expr=\coordindex, y=A3, col sep=comma] {\datatable};
\legend{A1: ISAM for Satellites, A2: Space Exp. Campaigns, A3: Satellite Constellations}

\end{axis}
\node[above, align=center] at (current bounding box.north) {\figtitle};

\end{tikzpicture}

%% file: tikz_figs/2010s.tex

\def\filename{2010s.csv}

\pgfplotstableread[col sep=comma]{2010s.csv}\datatable

\def\figtitle{2010-2020}
\begin{tikzpicture}
\begin{axis}[
    ybar stacked,
    bar width=.5cm,
    width=\textwidth,
    height=5cm,
    xtick=data,
    xticklabels from table={\datatable}{x_label},
    x tick label style={
        align=center, 
        },
    xtick style={draw=none}, 
    enlarge x limits=0.03, 
    ylabel={Number of citations},
    ymin=0,
    ymax=120,
    ytick={0,20,40,60,80,100,120},
    ymajorgrids, 
    ytick style={draw=none}, 
    legend style={
        at={(1,0.98)}, 
        anchor=north east, 
        legend columns=1, 
        draw=none, 
        },
    ]

\addplot+[] table [x expr=\coordindex, y=A1, col sep=comma] {\datatable};
\addplot+[] table [x expr=\coordindex, y=A2, col sep=comma] {\datatable};
\addplot+[] table [x expr=\coordindex, y=A3, col sep=comma] {\datatable};
\end{axis}
\node[above, align=center] at (current bounding box.north) {\figtitle};
\end{tikzpicture}

%% file: tikz_figs/2020s.tex

\def\filename{2020s.csv}

\pgfplotstableread[col sep=comma]{2020s.csv}\datatable

\def\figtitle{2020-2023}
\begin{tikzpicture}
\begin{axis}[
    ybar stacked,
    bar width=.5cm,
    width=\textwidth,
    height=5cm,
    xlabel={References and Used Methods},
    xlabel style={at={(axis description cs:0.45,-0.25)}}, 
    xtick=data,
    xticklabels from table={\datatable}{x_label},
    x tick label style={
        align=center, 
        },
    xtick style={draw=none}, 
    enlarge x limits=0.03, 
    ylabel={Number of citations},
    ymin=0,
    ymax=40,
    ytick={0,10,20,30,40},
    ymajorgrids, 
    ytick style={draw=none}, 
    legend style={
        at={(1,0.98)}, 
        anchor=north east, 
        legend columns=1, 
        draw=none, 
        },
    ]

\addplot+[] table [x expr=\coordindex, y=A1, col sep=comma] {\datatable};
\addplot+[] table [x expr=\coordindex, y=A2, col sep=comma] {\datatable};
\addplot+[] table [x expr=\coordindex, y=A3, col sep=comma] {\datatable};
\end{axis}
\node[above, align=center] at (current bounding box.north) {\figtitle};
\end{tikzpicture}